\documentclass[11pt]{amsart}
\usepackage{amscd, amsmath, amsthm, amssymb, xy, xypic, multirow}
\newtheorem{theorem}{Theorem}[section]
\newtheorem{lem}[theorem]{Lemma}
\newtheorem{prop}[theorem]{Proposition}
\newtheorem{conjecture}[theorem]{Conjecture}
\newtheorem{condition}[theorem]{Condition}
\newtheorem{cor}[theorem]{Corollary}
\theoremstyle{definition}     
\newtheorem{defi}[theorem]{Definition}

\newtheorem{example}[theorem]{Example}

\newtheorem{remark}[theorem]{Remark}
\numberwithin{equation}{section}

\def \hd #1 {\bfseries #1  \mdseries}
\def \italic #1 {\bfseries \it #1 \rm \mdseries}
\def \ra {\rightarrow}
\def \cen #1 { \begin{center} #1 \end{center}}

\def \mbz {\mathbb Z}

\def \mbc {\mathbb C}
\def \mbp {\mathbb P}

\def \mbq  {\mathbb {Q}}
\def \mco  {\mathcal {O}}

\def \mcv  {\mathcal {V}}
\def \mck  {\mathcal {K}}
\def \mcp  {\mathcal {P}}
\def \Q {${\mathbb {Q}}\,$}
\def \Max {{\rm{Max} }}
\def \Pic {{\rm{Pic}}}
\def \rk {{\rm{rk}}}

\def \ch {{\rm{ch}}}

\def \Min {{\rm{Min}}}

\begin{document}
\title[Calabi--Yau construction]
{Calabi--Yau construction by smoothing normal crossing varieties }

\author[N.-H. Lee]{Nam-Hoon Lee }
\address{ School of Mathematics, Korea Institute for Advanced Study, Dongdaemun-gu, Seoul 130-722, Korea }
\email{nhlee@kias.re.kr}
\begin{abstract}
We investigate a method of construction of Calabi--Yau manifolds, that is, by smoothing normal crossing varieties. We develop some theories for calculating the Picard groups of the Calabi--Yau manifolds obtained in this method. Some applications are included, such as construction of new examples of Calabi--Yau 3-folds with Picard number one with some interesting properties.
\end{abstract}
\maketitle
\setcounter{section}{-1}
\section{Introduction}
A \italic{Calabi--Yau manifold} is a compact K\"ahler manifold with
trivial canonical class such that the intermediate cohomology groups of
its structure sheaf are trivial  ($h^i (\mco_X) =0$ for $0 < i <
\dim (X)$). Simple examples are a smooth quintic hypersurface in
$\mbp^4$ and a complete intersection of two cubics in $\mbp^5$.
Actually the majority of  the known examples come in this ways --- complete
intersections in toric varieties.
Here we investigate a rather unfamiliar construction of Calabi--Yau manifolds, that is,   by smoothing normal crossing varieties.
 A \italic{normal crossing
variety} is a reduced complex analytic space which is locally
isomorphic to a normal crossing divisor on a smooth variety. This construction is intrinsically up to smooth deformation.
We regard two Calabi--Yau manifolds as  the same one  if they are connected by a chain of smooth deformations.

The starting point is the smoothing theorem (Theorem \ref{kana}) of Y.
Kawamata and Y. Namikawa. The theorem guarantees  a smoothing of some
normal crossing variety with certain conditions to a Calabi--Yau
manifold. Seeing that most of known examples of Calabi--Yau 3-folds
arise in the toric setting, this smoothing theorem may be useful for
constructing families of non-toric Calabi--Yau manifolds. Despite the
potentiality of the theorem, the properties of Calabi--Yau manifolds
obtained by the theorem were rather obscure. Firstly we demonstrate how to calculate
the Picard groups and Chern classes.
 More generally we show:
\smallskip  \newline
{\emph {In a semistable degeneration, one can construct the
second integral cohomology group of the generic fiber from those of components in the central normal crossing fiber
at least up to torsion  if every cycle in it is analytic (i.e. $h^{0,2}(X_t)=0$). \smallskip  \newline}

Later, we investigate the construction of some
Calabi--Yau manifolds by smoothing. In Section \ref{exam} , we apply our methods to the examples of Kawamata and Namikawa (\cite{KaNa}).
We show that each of classes (up to the  topological
Euler numbers) of the examples   actually
consists of multiple different families of Calabi--Yau 3-folds.
In Section \ref{cyfano}, we construct 7 new examples of Calabi--Yau 3-folds with Picard number one and show that they are connected to other Calabi--Yau 3-folds by (necessarily non-smooth) projective flat deformation.

To explain the method of calculation, we deduce the following corollary from Theorem \ref{mainthm}, Proposition \ref{chm} and Proposition \ref{6.4}.  Under a condition (Condition \ref{consur}) that can be checked for each individual case, we have the following:
\begin{cor}\label{qthm} Let $X_0 = Y_1 \cup Y_2$ be a normal crossing, which is the central fiber in a semistable degeneration of a Calabi--Yau $n$-fold $Z$ with $n \geq 3$. Then
\begin{enumerate}
\item The Picard group $\Pic (Z)$ of $Z$ is isomorphic to
$$\{  (l_1, l_2)  \in H^2 (Y_1, \mbz) \times H^2 (Y_1, \mbz) \big| l_1|_{D} = l_2|_{D} \,\, {\rm in} \,\,H^2(D, \mbz)\}
/\langle   (D, -D)\rangle  $$
up to torsion with the cup product preserved, where $D = Y_1 \cap Y_2$.
\item Furthermore if
\cen{$\omega_{X_0} = \mco_{X_0}$.}
Then the second Chern class  $c_2 (Z)$ of $Z$ corresponds to  $ (c_2
(Y_1), c_2 (Y_2))$.
\end{enumerate}
\end{cor}
The second assertion will be explained later. When $n=3$, it means that for any $l \in \Pic (Z)$, we have
$$l \cdot c_2 (Z) =  (l_1, l_2) \cdot  (c_2 (Y_1), c_2 (Y_2)) = l_1 \cdot c_2 (Y_1) + l_2 \cdot c_2 (Y_2),$$
where $ (l_1, l_2)$ is an element of $H^2 (Y_1, \mbz) \times H^2 (Y_1, \mbz)$ that corresponds to $l$ in the isomorphism in the above corollary.

 We will demonstrate a quick application of it in the following section.
The readers who are mostly interested in Calabi--Yau construction may go over the quick example in the next section and then jump to sections (\ref{kncy}, \ref{cyfano}) that are specialized for Calabi--Yau construction.

If every component of a normal crossing variety is smooth, we say that the
normal crossing variety is simple.
From now on, a normal crossing variety will mean a simple one unless stated otherwise.
 We assume that the
Clemens--Schmid  sequence is exact for given semistable
degeneration.

\medskip
{\it Acknowledgments}

\medskip
The paper is a part of the author's thesis work. The author would like to express his sincere thanks to his
teacher, Prof.\ Igor V. Dolgachev for his
guidance, including many suggestions and encouragement.

\section{A quick example}\label{quick}

As a simple demonstration how to use Corollary \ref{qthm},  we
construct the Picard group and the Chern class of a Calabi--Yau
3-fold obtained by smoothing.
Let us consider one of the examples that Y. Kawamata and Y. Namikawa
considered  in their paper (\cite{KaNa}). Let $Y_1 = \mbp^3$, $D$ be
a smooth quartic in $\mbp^3$ and $Y_2 $ be the blow-up of $\mbp^3$
along a smooth curve $c \in |\mco_D(8)|$.  Let $X_0 = Y_1 \cup_D
Y_2$, where `$\cup_D$' means pasting along $D$ (Note that $Y_1$ and
$Y_2$ contain  copies of $D$). The theorem of Y. Kawamata and Y.
Namikawa guarantees that $X_0$ is smoothable to a Calabi--Yau 3-fold
$Z$, that is, $X_0$ is the central fiber in a semistable degeneration of $Z$. Note that $\omega_{X_0}=\mco_{X_0}$.

It was quite mysterious what kind of Calabi--Yau  it is
although its topological Euler number was calculated ($e (Z) =
-296$). Now we can characterize it more. By Corollary \ref{qthm},
$$\Pic (Z) = \langle  (H, \pi^* (H))\rangle  $$
up to torsion, where $H$ is the ample generator of $H^2 (\mbp^3, \mbz)$ and $\pi: Y_2 \ra \mbp^3$ is the blow-up with the exceptional divisor $E$. So the Calabi--Yau has Picard number one. Let $\rho= (H, \pi^* (H))$. Then
$$\rho^3 = H^3+\pi^* (H)^3 = 1+1 = 2$$
and
$$\rho \cdot c_2 (Z) = H \cdot c_2 (Y_1) +  \pi^* (H) \cdot c_2 (Y_2) = 6+ 38 = 44.$$
Since $\rk\Pic(Z) = 1$, the divisor class $\rho$ is  ample on $Z$. According to a theorem on Calabi--Yau 3-folds (Corollary 1, \cite{GaPu}),  $8\rho$ is very ample and $Z$ has a projectively normal embedding into
$\mbp^N$, where
\begin{align*}
N=h^0(Z, \mco_Z(8\rho)) -1 &= \chi (\mco_Z (8\rho)) -1 \\
                    &= t_{6}({\rm Todd}(Z) \cdot \ch(8\rho))-1\\
                    &=  \frac{1}{6}(8\rho)^3 +  \frac{1}{12} c_2(Z) \cdot 8\rho -1 \\
                    &= 199.
\end{align*}
Has there been known any Calabi--Yau 3-fold with these invariants ( {$\Pic(Z) = \,\langle \rho\rangle $,} $\rho^3 = 2$,
$\rho\cdot c_2(Z) = 44$ and $e(Z) = -296$)? Indeed there is such a Calabi--Yau 3-fold -- a degree 8 hypersurface in $\mbp (1,1,1,1,4)$ (see, for example, Table 1 in \cite{EnSt}). Denote it by
$$X (8) \subset \mbp (1,1,1,1,4).$$
Let us try to construct its degeneration as above. $X (8)$ is linearly equivalent to $X (4) + X (4)$ as divisors.
So we can construct
a degeneration of $X (8)$ to $W_1 \cup W_2$, where $W_i$ is isomorphic to $X (4)$.
Note that $X (4)$ is a copy of $\mbp^3$. The total space of the degeneration is not smooth but it admits a small resolution, which makes the central fiber  a normal crossing of
$W'_1$ and
$W'_2$, where  $W'_1$ is a copy of $\mbp^3$ and $W'_2$  is isomorphic to $Y_2$. Indeed, we can
construct a degeneration of $X (8)$ whose central fiber is isomorphic to $X_0 =
Y_1 \cup Y_2$.

\section{Notations and preliminaries}

 Let $\pi: X \ra \bar \Delta$ be a proper map from a  smooth $(n+1)$-fold $X$ with
boundary onto a closed disk $\bar \Delta = \{ t \in \mbc \big| \|t\|
\leq 1\}$ such that the fibers $X_t = \pi^{-1} (t)$ are connected
K\"ahler  $n$-folds for every $t \neq 0$ (generic) and the central
fiber $X_0=\pi^{-1} (0) =\bigcup_\alpha Y_\alpha$ is a normal
crossing of  $n$-folds. We denote the generic fiber by $X_t$. The
condition, $t \neq 0$, is assumed in this notation. We call such a
map $\pi$ (or simply the total space $X$) a semi-stable degeneration
of $X_t$ and we say that $X_0$ is smoothable to $X_t$ with the
smooth total space. Let $[X_t]$ and $[Y_\alpha]$ be the fundamental
classes of $X_t$ and $Y_\alpha$ in $H^{2n} (X,
\partial X; \mbz)$ respectively. Then
$$[X_t] = [X_0] = \sum_\alpha [Y_\alpha].$$

Let $Y_{ij} = Y_i \cap Y_j$, $Y_{ijk} = Y_i \cap Y_j \cap Y_k$ and etc. Consider
the exact sequence
$$0 \ra \mbz_{X_0} \ra \bigoplus_\alpha \mbz_{Y_\alpha} \ra \bigoplus_{i <  j}
\mbz_{Y_{ij}} \stackrel{\tau}\ra  \bigoplus_{i <  j < k} \mbz_{Y_{ijk}} {\ra}
\cdots  $$
We introduce a short exact sequence,
$$0 \ra \mbz_{X_0} \ra \bigoplus_\alpha \mbz_{Y_\alpha} \ra  {{\rm{ker}}} (\tau)
{\ra} 0 $$
to give an exact sequence
\begin{align}
 \cdots \ra H^m (X_0, \mbz) \stackrel{\psi_m}{\ra}
\bigoplus_\alpha H^m  (Y_\alpha, \mbz) \ra
H^{m} ({{\rm{ker}}} (\tau))    \ra     \cdots   \label{seq}
\end{align}
Note that $\psi_m = \bigoplus_\alpha
j_\alpha^*|_{H^m (X_0, \mbz)}$, where $j_\alpha: Y_\alpha
\hookrightarrow X_0$ is the inclusion.

\begin{remark} \label{retract} $X_0$ is a deformation retract of $X$  (\cite{Cl}).
\end{remark}

So we have isomorphisms
$$k^*: H^m (X_0, \mbz) \widetilde \ra H^m (X, \mbz),$$
where $k: X_0 \hookrightarrow X$ is the inclusion. Let $i_\alpha: Y_\alpha
\hookrightarrow X$ be the inclusion.

Since $i_\alpha = k \circ j_\alpha$, we have the following commutative diagram:
$$
\xymatrix{
 &H^m (X, \mbz) \ar[dr]^{\varphi_m} && \\
 H^{m-1} ({{\rm{ker}}} (\tau)) \ar[r]& H^m (X_0, \mbz) \ar[u]^{k^*}_{\widetilde =}
\ar[r]^{\psi_m} &\bigoplus_\alpha H^m  (Y_\alpha, \mbz) \ar[r]&
H^{m} ({{\rm{ker}}} (\tau)),
}$$
where $\varphi_m = \bigoplus_\alpha i_\alpha^*|_{H^m (X, \mbz)}$. Let
$$G^m (X_0, \mbz) = {\rm{im}} (\varphi_m) = {\rm{im}} (\psi_m) =
{{\rm{ker}}}\left (\bigoplus_\alpha H^m
 (Y_\alpha, \mbz)\ra H^{m} ({{\rm{ker}}} (\tau))\right).$$
Let $\phi_m:H^m (X, \mbz) \ra G^m (X_0, \mbz)$ be the map which is derived from
$\varphi_m$.
By definition of $G^m (X_0, \mbz)$, $\phi_m$ is surjective.
Note that $G^m (X_0, \mbz)$ depends only on $X_0$ and makes sense for general
normal crossings as well.
We will construct $H^2 (X_t, \mbz)$ and $H^{2n-2}  (X_t, \mbz)$ as quotients of $G^2 (X_0, \mbz)$ and $G^{2n-2} (X_0, \mbz)$ respectively.

For an abelian
group $A$, we set
$$A_f = A / T,$$
where $T$ is the torsion part of $A$. We also denote $G^m (X_0, \mbz) \otimes_\mbz \mbq$ by
$G^m (X_0, \mbq)$ and etc.

\section{Multilinear maps and results from the Clemens--Schmid sequence}
Consider non-negative integers,  $q_1, \cdots, q_k$ such that
$$q_1 + \cdots +q_k = n.$$
Then there are multilinear maps, defined by the cup-product:
$$H^{2q_1} (X_t, \mbz) \times \cdots \times H^{2q_k} (X_t, \mbz) \ra \mbz,$$

$$H^{2q_1} (Y_\alpha, \mbz) \times \cdots \times H^{2q_k} (Y_\alpha, \mbz) \ra \mbz.$$
The latter one induces a multilinear map
$$\bigoplus_{\alpha_1} H^{2q_1} (Y_{\alpha_1}, \mbz) \times \cdots \times \bigoplus_{\alpha_k} H^{2q_k} (Y_{\alpha_k}, \mbz) \ra \mbz,$$
setting the mixed terms to be equal to zero. By restricting,  we can define a multilinear map
$$G^{2q_1} (X_0, \mbz) \times \cdots \times G^{2q_k} (X_0, \mbz) \ra \mbz.$$
Furthermore we define a multilinear map
$$H^{2q_1} (X, \mbz) \times \cdots \times H^{2q_k} (X, \mbz) \ra \mbz$$
by setting
$$ (l_1, \cdots, l_k) \mapsto l_1 \cdot ... \cdot l_k \cdot [X_t].$$

\begin{prop}\label{forms} These maps are compatible with the following diagram:

$$
\xy
 (0,0)*{H^{2q_1} (X, \mbz) \times \cdots \times H^{2q_k} (X, \mbz)}="1";
 (-35,-20)*{H^{2q_1} (X_t, \mbz) \times \cdots \times H^{2q_k} (X_t, \mbz)}="2";
 (35,-20)*{ G^{2q_1} (X_0, \mbz) \times \cdots \times G^{2q_k} (X_0, \mbz), }="3";

{\ar "1";"2"}
{\ar^{\phi_{2q_1} \times \cdots \times \phi_{2q_k}} "1";"3"}
\endxy
$$
where the map
$$H^{2q_1} (X, \mbz) \times \cdots \times H^{2q_k} (X, \mbz) \ra
H^{2q_1} (X_t, \mbz) \times \cdots \times H^{2q_k} (X_t, \mbz)$$
is induced by the inclusion, $i:X_t \hookrightarrow X$.
\end{prop}
\begin{proof}
For $ (l_1, \cdots, l_k) \in H^{2q_1} (X, \mbz) \times \cdots \times H^{2q_k} (X, \mbz)$,
\begin{align*}
i^* (l_1) \cdot ... \cdot i^* (l_k) &= \int_{X_t}i^* (\omega_{l_1}) \wedge \cdots \wedge i^* (\omega_{l_k}) \\
                         &=\int_{X}\omega_{[X_t]} \wedge \omega_{l_1} \wedge \cdots \wedge \omega_{l_k} \\
                         &=l_1 \cdot ... \cdot l_k \cdot [X_t],
\end{align*}
where $\omega_{l_k}$ and $\omega_{[X_t]}$ are the De Rham classes that correspond to $l_k \in H^{2q_k} (X, \mbz)$ and $[X_t] \in H^2 (X, \partial X; \mbz)$ respectively. Similarly,
\begin{align*}
\phi_{q_1}  (l_1) \cdot ... \cdot \phi_{q_k}  (l_k) &= \sum_\alpha \int_{Y_\alpha} i_\alpha^* (\omega_{l_1}) \wedge \cdots \wedge i_\alpha^* (\omega_{l_k})\\
                                       &=\sum_\alpha \int_{X}\omega_{[Y_\alpha]} \wedge \omega_{l_1} \wedge \cdots \wedge \omega_{l_k}\\
                                       &=\sum_\alpha l_1 \cdot ... \cdot l_k \cdot [Y_\alpha] \\
                                       &=l_1 \cdot ... \cdot l_k \cdot \sum_\alpha[Y_\alpha] \\
                                       &=l_1 \cdot ... \cdot l_k \cdot [X_0]\\
                                       &=l_1 \cdot ... \cdot l_k \cdot [X_t].
\end{align*}
\end{proof}

\begin{remark}\label{formsq} The above proposition also holds with rational coefficients.
\end{remark}

By analyzing the Clemens--Schmid exact sequence, we have
\begin{theorem}\label{css} Suppose that every cycle in $H^2 (X_t, \mbz)$ is analytic, i.e. $h^{2,0} (X_t) =0$. Then the maps
$$H^2 (X, \mbq) \ra H^2 (X_t, \mbq),$$
$$H^{2n-2} (X, \mbq) \ra H^{2n-2} (X_t, \mbq)$$
are surjective and
$$h^2 (X_t) = h^2 (X_0) - r +1,$$
where $r$ is the number of components of $X_0$.
\end{theorem}
For detailed proof, we refer to Chapter III of \cite{Lee}.

\section{Construction of $H^{2} (X_t, \mbz)$}

We construct $H^{2} (X_t, \mbz)$ under the assumption that the  following map
$$H^2 (X, \mbz) \ra H^2 (X_t, \mbz)$$
is surjective up to torsion.

Let $[Y_\alpha]'$ be the image of $[Y_\alpha]$ under the map
$$ H^2 (X, \partial X ; \mbz) \ra H^2 (X, \mbz).$$
Note that $[X_t]$ goes to zero by this map, i.e. $[X_t]'=0$.

Now we introduce a degenerate subgroup $N^2(X_0, \mbz)$ of $G^2(X_0,
\mbz)$, which is the main point. We will obtain $H^2(X_0, \mbz)$ by taking a quotient of
$G^2(X_0, \mbz)$ over $N^2(X_0, \mbz)$. Note the canonicity of its
form.
\begin{lem} \label{rank} Let $X_0 = Y_1 \cup \cdots \cup Y_r$ and  $e_i$ be the image of $[Y_i]'$ under the map
$$\phi_2: H^2 (X, \mbz) \ra G^2 (X_0, \mbz).$$ Let $NG^2(X_0, \mbz)$ be a subgroup of $G^2 (X_0, \mbz)$, generated by $e_1, \cdots, e_r$. Then then rank of $NG^2(X_0, \mbz)$ is $r-1$.
\end{lem}
\begin{proof} Note that $e_i$ is the image of $[Y_i]$ under the map
$$H^2 (X, \partial X ; \mbz) \ra G^2 (X_0, \mbz),$$
which is induced by
$$H^2 (X, \partial X ; \mbz) \ra H^2 (Y_\alpha, \mbz).$$
Since
$$e_1 + \cdots + e_r = \phi_2 ([Y_1]' + \cdots + [Y_r]') = \phi_2 ([X_0]') = \phi_2 ([X_t]') = 0,$$
the rank of $NG^2(X_0, \mbz)$ is less than $r$. Suppose that the rank is less than {\mbox{$r-1$,}} then there are integers
$a_1, \cdots, a_r$ such that  $ (a_1, \cdots, a_r)$ is not a multiple of $ (1, \cdots, 1)$ and that
$$a_1 e_1 + \cdots + a_r e_r = 0.$$
Since $e_1 + \cdots + e_r = 0$, we can add or subtract $ (e_1 + \cdots + e_r)$ from $a_1 e_1 + \cdots + a_r e_r$ as many times as we want. So we may assume that all the $a_j$'s  are non-negative and $a_i = 0$ for some $i$.
The projection of $a_1 e_1 + \cdots + a_r e_r$ to $H^2 (Y_i, \mbz)$ is the first Chern class of an effective divisor,
$$\sum_{j (\neq i)} a_j Y_j \cap Y_i$$
on $Y_i$.  Accordingly we have
\cen{$a_j = 0$ or $Y_j \cap Y_i = \emptyset$.}
So if $Y_j \cap Y_i \neq \emptyset$, then $a_j = 0$. Taking another projection to $H^2 (Y_j, \mbz)$, we can repeat the same argument to conclude that if $Y_k$ belongs to a connected component of $X_0$ that contains $Y_i$, then $a_k = 0$. But $X_0$ is connected because $X_t$ is so. Therefore we have
\cen{$a_k = 0$ for any $k$,}
which is a contradiction. So the rank should be equal to $r-1$.
\end{proof}

By the above lemma,
$$NG^2(X_0) = \langle  e_1, \cdots, e_{r-1}\rangle .$$
From
$$0=[X_t]|_{Y_i} = [X_0]_{Y_i} = [Y_1] |_{Y_i} + \cdots + [Y_{i}]|_{Y_{i}} + \cdots +[Y_r] |_{Y_i},$$
we have
$$[Y_{i}]|_{Y_{i}} = - \sum_{j (\neq i)} [Y_j] |_{Y_i}.$$
So
$$e_i =  (Y_1 \cap Y_i, \cdots, Y_{i-1} \cap Y_{i}, - \sum_{j (\neq i)} Y_j \cap Y_i, Y_{i+1} \cap Y_i, \cdots,
Y_r \cap Y_i) \in \bigoplus_\alpha H^2 (Y_\alpha, \mbz).$$

\begin{lem} \label{nonzero}Assume that $h^{2,0} (X_t) = 0$ and that an element $l \in H^2 (X, \mbz)_f$ goes to a non-zero element $l' \in H^2 (X_t, \mbz)_f$. Then $l$ goes to a non-zero element $l''$  under the map $$H^2 (X, \mbz)_f \ra G^2 (X_0, \mbz)_f.$$
\end{lem}
\begin{proof} Since $l' \neq 0$, by Poincar\'e duality, there is an element $k' \in H^{2n-2} (X_t, \mbz)_f \subset H^{2n-2} (X_t, \mbq)$ such that $l' \cdot k' \neq 0$.  According to Theorem \ref{css}, the map
  $$H^{2n-2} (X, \mbq) \ra H^{2n-2} (X_t, \mbq)$$
  is surjective. So there is an element $k \in H^{2n-2} (X, \mbq)$ that goes to $k'$. Let $k''$ be the image of $k$ under the map
$$H^2 (X, \mbq) \ra G^2 (X_0, \mbq).$$
By Proposition \ref{forms} and Remark \ref{formsq}, we have
$$0 \neq l' \cdot k' = l \cdot k \cdot [X_t]= l'' \cdot k''.$$
So $l''$ cannot be zero.
\end{proof}

\begin{remark}\label{gzero2} Since $Y_i \cap X_t =\emptyset$, The image of $[Y_i]'$ under the map
$$H^2 (X, \mbz) \ra H^2 (X_t, \mbz)$$
is zero.
\end{remark}

\begin{prop} \label{main} Suppose that $H^2 (X, \mbz) \ra H^2 (X_t, \mbz)$ is surjective up to torsion and that $h^{2,0} (X_t) = 0$. Then
$H^2 (X, \mbz)$ and $H^2 (X_t, \mbz)$ are isomorphic to $G^2 (X_0, \mbz)$ and $G^2 (X_0, \mbz) / NG^2(X_0, \mbz)$ respectively up to torsion with the cup products preserved.
\end{prop}
\begin{proof}  By Lemma \ref{nonzero}, Remark \ref{gzero2} and the surjectivity of the map
$$H^2 (X, \mbz) \ra G^2 (X_0, \mbz),$$
we have
$$\rk (H^2 (X_t, \mbz)) \leq \rk (G^2 (X_0, \mbz) / NG^2(X_0, \mbz)).$$
By Theorem \ref{css},
$$\rk (H^2 (X_t, \mbz)) = \rk (H^2 (X_0, \mbz)) -r +1.$$
Therefore,
\begin{align*}
\rk (H^2 (X_t, \mbz)) &\leq \rk (G^2 (X_0, \mbz) / NG^2(X_0, \mbz)) \\
                &= \rk (G^2 (X_0, \mbz)) -\rk (NG^2(X_0, \mbz)) \\
                &= \rk (G^2 (X_0, \mbz)) -r+1 \,\,\, (\because \,\, \rm{Lemma \,\, \ref{rank}})\\
                &\leq  \rk (H^2 (X, \mbz)) -r +1 \\
                &=  \rk (H^2 (X_0, \mbz)) -r +1 \,\,\, (\because \,\, \rm{Remark \,\, \ref{retract}})\\
                &= \rk (H^2 (X_t, \mbz)).
\end{align*}
So we have
$$\rk (H^2 (X, \mbz)) = \rk (G^2 (X_0, \mbz))$$
and
$$\rk (H^2 (X_t, \mbz)) = \rk (G^2 (X_0, \mbz) / NG^2(X_0, \mbz)).$$
The first equation implies that the map
$$H^2 (X, \mbz)\ra G^2 (X_0, \mbz)$$
is an isomorphism up to torsion.
Let $NH^2 (X, \mbz)$ be a subgroup of $H^2 (X, \mbz)$, generated by $[Y_1]', \cdots, [Y_r]'$. By Remark \ref{gzero2} and the surjectivity assumption, we have a surjective map
\begin{equation} \label{smap}
H^2 (X, \mbz)/NH^2 (X, \mbz) \ra H^2 (X_t, \mbz)_f.
\end{equation}
Since
\begin{align*}
\rk (H^2 (X, \mbz)/NH^2 (X, \mbz)) &= \rk (H^2 (X, \mbz)) - \rk (NH^2 (X, \mbz))\\
                           &= \rk (G^2 (X_0, \mbz)) - \rk (NG^2(X_0, \mbz))\\
                           &=\rk (H^2 (X_t, \mbz)),
\end{align*}
map (\ref{smap}) is an isomorphism up to torsion.   Therefore
$$\left (G^2 (X_0, \mbz) / NG^2(X_0, \mbz)\right)_f \simeq \left (H^2 (X, \mbz)/NH^2 (X, \mbz)\right)_f \simeq
H^2 (X_t,\mbz)_f.$$
By Proposition \ref{forms}, the cup products are preserved in any of maps in this proof.

\end{proof}

We denote $G^2 (X_0, \mbz)/ NG^2(X_0, \mbz)$ by $RG^2 (X_0, \mbz)$.

\section{Construction of $H^{2n-2} (X_t, \mbz)$ and circumventing the assumption of surjectivity}

In the previous section, we have constructed $H^{2} (X_t, \mbz)$.
 Now we construct $H^{2n-2} (X_t, \mbz)$ as a quotient of
$G^{2n-2} (X_0, \mbz)$.

Let $ NH^{2n-2} (X, \mbz) = \{ l \in H^{2n-2} (X, \mbz) \big| l \cdot l' \cdot [X_t] =0,
\forall l' \in H^{2} (X, \mbz) \}$.
\begin{lem}\label{p5.2} Suppose that $h^{2, 0} (X_t) =0$, then
$${{\rm{ker}}} (  H^{2n-2} (X, \mbz) \ra H^{2n-2} (X_t, \mbz)_f) = NH^{2n-2} (X, \mbz).$$
\end{lem}
\begin{proof}  Let
$$l \in \ker (H^{2n-2} (X, \mbz) \ra H^{2n-2} (X_t, \mbz)_f).$$
For any $l' \in H^{2} (X, \mbz)$, we have
$$l \cdot l' \cdot [X_t] = i^* (l) \cdot i^* (l')  =0$$
by Proposition \ref{forms}. So $l \in NH^{2n-2} (X, \mbz)$.

 On the other hand, let $l \in  NH^{2n-2} (X, \mbz)$.   Since
$$H^2 (X, \mbq) \ra H^2 (X_t, \mbq)$$ is surjective  (Theorem \ref{css}), for any $\alpha \in H^{2} (X_t, \mbq)$, there is an element $l' \in H^2 (X, \mbq)$ that goes to   $\alpha$. Then

 $$i^* (l) \cdot \alpha = l \cdot l' \cdot [X_t] = 0.$$
 Since Poincar\'e pairing on
 $$H^{2n-2} (X_t, \mbq) \times H^{2} (X_t,\mbq)$$
 is non-degenerate, $i^* (l) = 0$ in  $H^{2n-2} (X_t, \mbq)$. But this means that $i^* (l) = 0$ in $H^{2n-2} (X_t, \mbz)_f$, i.e
$$l \in {{\rm{ker}}} ( H^{2n-2} (X, \mbz)_f \ra
H^{2n-2} (X_t, \mbz))_f.$$
\end{proof}

Let
$$NG^{2n-2} (X_0, \mbz) = \{ l \in G^{2n-2} (X_0, \mbz) \big| l \cdot l' =0, \forall l'
\in G^{2} (X_0, \mbz)\}.$$

\begin{prop}\label{5.4} Suppose that $h^{2,0} (X_t) = 0$ and
$$H^{2n-2} (X, \mbz)_f \ra H^{2n-2} (X_t, \mbz)_f$$
is surjective, then
$$H^{2n-2} (X_t, \mbz)_f \simeq G^{2n-2} (X_0, \mbz)/ NG^{2n-2} (X_0, \mbz).$$
\end{prop}
\begin {proof}
From Lemma \ref{p5.2} and the assumption of surjectivity, we have
$$H^{2n-2} (X_t, \mbz)_f \simeq H^{2n-2} (X, \mbz) / NH^2 (X, \mbz).$$
Note that the map
$$\phi_{2n-2}: H^{2n-2} (X, \mbz) \ra G^{2n-2} (X_0, \mbz)$$
is surjective. So it is enough to   show that
$$ NH^{2n-2} (X, \mbz) = \phi_{2n-2}^{-1} (NG^{2n-2} (X_0, \mbz)).$$
Let $l \in NH^{2n-2} (X, \mbz)$.  For any $m \in
G^{2} (X_0,
\mbz)$, there is  an element $l' \in
H^{2} (X, \mbz)$ such that $ \phi_{2} (l') = m$. Then
$$\phi_{2n-2} (l) \cdot m = l \cdot l' \cdot [X_t] = 0.$$
So
$$\phi_{2n-2} (l) \in NG^{2n-2} (X_0, \mbz),$$
i.e.
$$ l \in  \phi_{2n-2}^{-1} (NG^{2N-2} (X_0, \mbz)).$$
We conclude that
$$ NH^{2n-2} (X, \mbz) \subset \phi_{2n-2}^{-1} (NG^{2n-2} (X_0, \mbz)).$$
Conversely, let $l \in \phi_{2n-2}^{-1} (NG^{2n-2} (X_0, \mbz))$, i.e.
$$\phi_{2n-2} (l) \in NG^{2n-2} (X_0, \mbz).$$
Suppose $ l \notin  NH^{2n-2} (X, \mbz)$, then there is an element $l' \in H^{2} (X, \mbz)$
such
that $l \cdot l' \cdot [X_t] \neq 0$.   But then
$$0 = \phi_{2n-2} (l) \cdot \phi_{2n-2} (l') = l \cdot l' \cdot [X_t] \neq 0,$$
which is a contradiction. Therefore $l \in  NH^{2n-2} (X, \mbz)$ and
$$ NH^{2n-2} (X, \mbz) \supset \phi_{2n-2}^{-1} (NG^{2n-2} (X_0, \mbz)).$$
\end{proof}

We denote $G^{2n-2} (X_0, \mbz)/ NG^{2n-2} (X_0, \mbz)$ by
$RG^{2n-2} (X_0, \mbz)$.

In Proposition \ref{main} and Proposition \ref{5.4}, we assumed the surjectivity of the following maps:
\begin{equation} \label{map2}
H^2  (X, \mbz)_f \ra  H^2  (X_t, \mbz)_f \,\,\, {\rm and} \,\,\,
H^{2n-2}  (X, \mbz)_f \ra  H^{2n-2}  (X_t, \mbz)_f.
\end{equation}

This assumption is supported by the surjectivities of the maps:
\begin{equation} \label{map2q}
H^2  (X, \mbq) \ra  H^2  (X_t, \mbq) \,\,\, {\rm and} \,\,\,
H^{2n-2}  (X, \mbq) \ra  H^{2n-2}  (X_t, \mbq)
\end{equation}
in Theorem \ref{css} when $h^{2,0} (X_t) =0$.

The following proposition shows how to deal with this assumption and how to produce
non-conjectural result in each case.

\begin{prop} \label{check}When $h^{2,0} (X_t) =0$, the maps (\ref{map2})  are surjective
if and only of
 the pairing on $$RG^2  (X_0, \mbz) \times RG^{2n-2}  (X_0, \mbz)$$
is unimodular.
\end{prop}
\begin{proof}
$\Rightarrow$: If the maps (\ref{map2}) are surjective, then

\cen{$RG^2  (X_0, \mbz)_f \simeq H^2 (X_t, \mbz)_f $ and
$RG^{2n-2}  (X_0, \mbz) \simeq H^{2n-2} (X_t, \mbz)_f $}
by Proposition \ref{main} and Proposition \ref{5.4}.
The unimodularity follows from that of the pairing on
$$H^2  (X_t, \mbz) \times H^{2n-2}  (X_t, \mbz).$$
$\Leftarrow$: Suppose that the pairing on $$RG^2  (X_0, \mbz) \times RG^{2n-2}  (X_0, \mbz)$$
is unimodular.
  By Theorem \ref{css}, the maps (\ref{map2q}) are surjective. So  the cokernels of  the maps (\ref{map2}) have
no free part.  By  Poincar\'e duality, the pairing on $$H^2  (X_t, \mbz)_f
\times H^{2n-2}  (X_t, \mbz)_f$$
is unimodular. Suppose that either of the cokernels has a non-zero element,
which must be torsion. Then the pairing on $$RG^2  (X_0, \mbz) \times
RG^{2n-2}  (X_0, \mbz)$$
cannot be unimodular, which contradicts our assumption. So
the cokernels are all trivial, i.e. the maps (\ref{map2}) are surjective.
\end{proof}
So we have a checkable equivalent condition for the surjectivity.

\begin{condition}\label{consur}
The pairing on $$RG^2  (X_0, \mbz) \times RG^{2n-2}  (X_0, \mbz)$$
is unimodular.
\end{condition}

Now we state the main theorem of this paper, which follows immediately from Proposition \ref{main}, Proposition {5.4} and Proposition \ref{check}

\begin{theorem} \label{mainthm} Suppose that  $h^{2,0}(X_t) =0$ and  that Condition \ref{consur} is satisfied. Then
$H^2(X_t, \mbz)$ and $H^{2n-2}(X_t, \mbz)$ are isomorphic to $RG^2(X_0, \mbz)$ and $RG^{2n-2}(X_0, \mbz)$  up to torsion respectively.
\end{theorem}

Condition \ref{consur} will be automatically assumed in the rest of this paper and we will not state it explicitly.
All the examples the author considered satisfy Condition \ref {consur} when $h^{2,0} (X_t)=0$. We make
a conjecture.

\begin{conjecture}When $h^{2,0} (X_t) = 0$, Condition \ref{consur} holds or equivalently the maps
$$H^2  (X, \mbz)_f \ra  H^2  (X_t, \mbz)_f$$
and
$$H^{2n-2}  (X, \mbz)_f \ra  H^{2n-2}  (X_t, \mbz)_f$$
are surjective
\end{conjecture}

With assumption of surjectivity of the maps
$$H^{2m} (X, \mbz)_f \ra H^{2m} (X_t, \mbz)_f,$$
 one can construct $H^{2m} (X_t, \mbz)$ up to torsion for general $m$ in a similar way.   We do not study it here  because the surjectivity is not obtained very often unless $m = 1$ or $n-1$.
 We refer interested readers to Chapter IV of \cite{Lee}.

Later, we will see many examples where $X_0$ is composed of two components

\begin{cor} \label{two1} If $X_0 = Y_1 \cup Y_2$ and $h^{2,0}(X_t) = 0$. Then
  $$NG^2(X_0, \mbz) =\langle e_1\rangle  = \langle  (D, -D)\rangle $$
  and $H^2 (X_t, \mbz)$ is isomorphic to
  $$\{ (l_1, l_2) \in H^2 (Y_1, \mbz) \times H^2 (Y_2, \mbz)
  \big| l_1|_{D} = l_1|_{D}\,\, {\rm{in}} \,\, H^2 (D, \mbz) \}/\langle D, -D\rangle $$
  up to torsion, where $D = Y_1 \cap Y_2$.
\end{cor}
\begin{proof} It is easy to see that
$$G^2 (X_0, \mbz) =\{ (l_1, l_2) \in H^2 (Y_1, \mbz) \times H^2 (Y_2, \mbz)
  \big| l_1|_{D} = l_1|_{D}\,\, {\rm{in}} \,\, H^2 (D, \mbz) \}$$
from the exact sequence  (\ref{seq}),
$$\cdots \ra H^2 (X_0, \mbz) \ra H^2 (Y_1, \mbz) \oplus H^2 (Y_2, \mbz) \ra H^2 (D, \mbz) \ra \cdots.$$

\end{proof}

\section{Chern Classes}

 We locate the Chern class of $X_t$ as an element of $\bigoplus_i G^{2i}  (X_0, \mbz)$  (what this statement means will be clear in Proposition \ref{chm}).

Note that the normal bundle $N_{X_t}$ of $X_t$ in $X$ is trivial.
Let $[Y_{\beta\alpha}^\alpha] = i_\alpha^* ([Y_\beta]) \in H^2 (Y_\alpha, \mbz)$.
By
the adjunction
formula,
$$i^*  (c  (X)) =   (1 + i^*  ([X_t])) \cdot c  (X_t) = c  (X_t).$$
So the Chern classes of $X_t$ come from the total space $X$.
$$i_\alpha^*  (c  (X)) =   (1+i_\alpha^*  ([Y_\alpha])) \cdot c  (Y_\alpha).$$
Since $X_t \cap Y_\alpha = \emptyset$,
\begin{align*}
0=i_\alpha^*  ([X_t]) &= i_\alpha^*  ([X_0])\\
                  &= i_\alpha^*  (\sum_\beta [Y_\beta]) \\
                  &= i_\alpha^*  ([Y_\alpha]) + \sum_{\beta   (\neq \alpha)}
[Y_{\beta\alpha}^\alpha].
\end {align*}
So
$$i_\alpha^*  ([Y_\alpha])  =- \sum_{\beta   (\neq \alpha)}
[Y_{\beta\alpha}^\alpha].$$
Let us locate $c  (X_t)$. Since
$$ c  (X_t) = i^*  (c  (X)),$$
$c  (X_t)$ corresponds to $\phi  (c  (X))$ in  $\bigoplus_i G^{2i}  (X_0, \mbz)$.
Let us calculate $\phi  (c  (X))$
\begin{align*}
\phi  (c  (X)) &= \sum_\alpha i_\alpha^*  (c  (X))\\
            &= \sum_\alpha   (1^{  (\alpha)}+i_\alpha^*  ([Y_\alpha]) \cdot c  (Y_\alpha))
\\
            &=\sum_{\alpha} \left   (1^{  (\alpha)}-\sum_{\beta  (\neq
\alpha)}[Y_{\beta
\alpha}^\alpha] \right )\cdot c  (Y_\alpha),
\end{align*}
where $1^{  (\alpha)}$ is the generator of $H^0  (Y_\alpha, \mbz)$.
$\phi  (c  (X))$ depends only on $X_0$. So it makes sense for general normal crossing
as well as the central fiber of a degeneration. Let us give a
definition
of the Chern class of a general normal crossing.
\begin{defi} Let $Z = \bigcup_\alpha W_\alpha$ be a normal crossing. Then the
Chern class of $Z$ is defined as an element of  $\bigoplus_{p,
\alpha} H^{2p}  (W_\alpha, \mbz)$ by
$$c  (Z) =  \sum_{\alpha} \left   (1^{  (\alpha)}-\sum_{\beta  (\neq
\alpha)}[W_{\beta\alpha}^\alpha] \right )\cdot c  (W_\alpha).$$
\end{defi}
Note that $c  (Z) \not \in \bigoplus_p  G^{2p}  (Z, \mbz)$ in general unless
$Z$ is the central fiber of a semistable degeneration.

For a smooth variety, the first Chern class is the anticanonical class of the
variety. Let us check what happens in the normal crossing case.
$$c_1  (Z) = - \sum_{\alpha \not = \beta} [W_{\beta\alpha}^\alpha]  +
\sum_\alpha  c_1  (W_\alpha) = - \sum_\alpha \omega_{Z}|_{W_\alpha},$$
where $\omega_{Z}$ is the canonical class for $Z$. So  the above definition
generalizes the conventional Chern class at least for the first Chern class.

Now we come back to our original topic. The following proposition, following from Proposition \ref{forms} and the definition of $c (X_0)$, provides us with a link between $c (X_t)$ and $c (X_0)$.
\begin{prop} \label{chm} Consider non-negative integers,  $q_1, \cdots, q_k, p$ and $i$ such that
$$q_1 + \cdots +q_k+ip = n.$$
Let $l_j \in H^{2 q_j } (X, \mbz)$ for $j = 1, \cdots, k$. Then we
have the following equality between two integers:
$$i^* (l_1) \cdot ... \cdot i^* (l_k) \cdot c_i (X_t)^p = \phi (l_1)\cdot ... \cdot \phi (l_k) \cdot c_i (X_0)^p.$$
\end{prop}

When $n=3$, we can construct all the even-dimensional cohomology groups of the generic fiber together with the Chern class from the geometry of the central fiber. Let $RG^0(X_0, \mbz) = G^0(X_0, \mbz)$. Let
$$NG^{2n} (X_0, \mbz) = \{ l \in G^{2n} (X_0, \mbz) \big| l \cdot l' =0, \forall l'
\in G^{0} (X_0, \mbz)\}$$
and $RG^{2n}(X_0, \mbz) = G^{2n} (X_0, \mbz)/NG^{2n} (X_0, \mbz)$.
 The following theorem follows immediately from what we have shown so far.

\begin{theorem} Suppose that $h^{2,0}(X_t) = 0$ and $n=3$, then we have an isomorphism:
\begin{align*}
&\left (H^0(X_t, \mbz) \oplus H^2(X_t, \mbz) \oplus H^4(X_t, \mbz) \oplus H^6(X_t, \mbz)  , c(X_t)\right )
\simeq \\
&\left (RG^0(X_0, \mbz) \oplus RG^2(X_0, \mbz) \oplus RG^4(X_0, \mbz) \oplus RG^6(X_0, \mbz)  , c(X_0)\right )
\end{align*}
up to torsion with the cup product preserved.
\end{theorem}

The following result for $ c_2  (X_0)$ is very useful in practice.

\begin{prop}\label{6.4} Suppose that the central fiber $X_0$ is composed of two components $Y_1$ and $Y_2$ and that  $\omega_{X_0}^{\otimes k} = \mco_{X_0}$ for some positive integer $k$.
Then
 $$c_2  (X_0) \simeq  c_2  (Y_1) + c_2 (Y_2)$$
 modulo an element that is trivial w.r.t.\ cup product.
\end{prop}
\begin{proof}
The condition
 $\omega_{X_0}^{\otimes k} = \mco_{X_0}$ implies that $c_1 (Y_1) \simeq [Y_{2\, 1}^1] $ and
$c_1 (Y_2) \simeq [Y_{1\, 2}^2] $ up to torsion.
So we have
\begin{align*}
  c_2  (X_0) &= c_2  (Y_1) +  c_2  (Y_1)   - [Y_{2\, 1}^1] \cdot c_1  (Y_1) -[Y_{1\, 2}^2] \cdot c_1  (Y_2) \\
               &\simeq c_2  (Y_1) +  c_2  (Y_1)   - [Y_{2\, 1}^1] \cdot [Y_{2\, 1}^1] -[Y_{1\, 2}^2] \cdot [Y_{1\, 2}^2]\\
               &\simeq c_2  (Y_1) +  c_2  (Y_1)   -  ([Y_{2\, 1}^1] - [Y_{1\, 2}^2])^2 \,\,\,\,\,\,\,\   (\because [Y_{2\, 1}^1] \cdot [Y_{1\, 2}^2]  =  0 )\\
               &\simeq c_2  (Y_1) +  c_2  (Y_1)  \,\,\,\,\,\,\,\   (\because [Y_{2\,1}^1 ]- [Y_{1\, 2}^2 ] \in NG^2 (X_0, \mbz).
\,\,\, {\rm{Corollary \,\,\ref{two1}}}).
\end{align*}
\end{proof}

\section{Calabi--Yau construction by smoothing normal crossing varieties}\label{kncy}
\label{exam} Now we state the theorem of Y. Kawamata and Y. Namikawa
(\cite{KaNa}).

\begin{theorem}[Y. Kawamata, Y. Namikawa] \label{kana}Let ${X_0} =
\bigcup_k Y_k$ be (not necessarily simple) normal crossing of dimension $n \geq 3$ such that
\begin{enumerate}
\item Its dualizing sheaf is trivial: $\omega_{X_0} = \mco_{X_0}$.
\item $H^{n-2}  (Y_k, \mco_{Y_k}) = 0$ for any $k$ and $H^{n-1}  ({X_0}, \mco_{X_0}) = 0$.
\item It is K\"ahler.
\item It has a logarithmic structure.
\end{enumerate}
Then ${X_0}$ is smoothable to a  $n$-fold with the smooth total space
(semistable degeneration).
\end{theorem}
The `logarithmic
structure' is also called as d-semistability. It will be explained shortly for a simpler case that is necessary in the paper.
For a Calabi--Yau $n$-fold $Z$ with $n \geq  3$, $\Pic (Z) \simeq H^2 (Z,
\mbz)$. So if the variety obtained by smoothing is a Calabi--Yau manifold, we can
construct its Picard group  from the central normal crossing. From now on, we restrict ourselves to the three dimension but our theory is readily generalized to higher dimension with minor modifications.
The case in which the central fiber has only two components occurs very
often.  That is, $X_0=Y_1 \cup Y_2$. Let $D = Y_1\cap Y_2$. Then the conditions in the above theorem correspond to the followings:
\begin{enumerate}
\item $D \in |{-}K_{Y_i}|$ for $i=1,2$.
\item $H^1(\mco_{Y_i}) = H^1(\mco_D)=0$. Note that this condition, together with (1),  implies that $D$ is a $K3$ surface.
\item There are ample divisors $H_1$, $H_2$ on $Y_1$, $Y_2$ respectively such that $H_1|_D$ is linearly equivalent to $H_2|_D$.
\item $N_{D/Y_1}\otimes N_{D/Y_2} = \mco_D$, where $N_{D/Y_i}$ is the normal bundle to $D$ in $Y_i$.  This condition is called d-semistability.
\end{enumerate}
By the above theorem, $X_0$ is smoothable to a  3-fold, $X_t$ with $K_{X_t} = 0$ and $h^{1}(\mco_{X_t}) = h^2(\mco_{X_t})=0$. Accordingly $X_t$ is a Calabi--Yau 3-fold.

Let us gather some facts in the following, which is a corollary of Theorem \ref{css}, Corollary \ref{two1} and Proposition \ref{6.4}:

\begin{cor}\label{9.1}
Let $X_t$ be the smoothing of $X_0$ as the above. Then
\begin{enumerate}
\item $h^{1, 1} (X_t) = h^{2} (Y_1) + h^{2} (Y_2) - k -1$, where
$$k = {\rm{rk}} ({\rm{im}} (H^2 (Y_1, \mbz) \oplus H^2 (Y_2, \mbz) \ra H^2 (D, \mbz))).$$
\item $h^{1, 2} (X_t) = 21 + h^{1, 2} (Y_1) + h^{1, 2} (Y_2) - k$.
\item ${\rm{Pic}} (X_t)$ is isomorphic to
$$\{ (l, l') \in \Pic (Y_1) \times \Pic (Y_2) \big|
l|_D = l'|_D \,\,{\rm in}\,\, \Pic (D)\} / \langle  (D, -D)\rangle $$
up to torsion.
\item $c_2 (X_0) = c_2 (Y_1) + c_2 (Y_2)$ modulo an element that is trivial w.r.t.\ cup product.
\end{enumerate}
\end{cor}
\begin{proof}\
\begin{enumerate}
\item From Theorem \ref{css},
\begin{align*} h^{1,1} (X_t)&= h^2 (X_t) \\
                    &= h^2 (X_0) - 2+1 \\
                     &=  (h^2 (Y_1) + h^2 (Y_2) - k) -1 \\
                     &=  h^{2} (Y_1) + h^{2} (Y_2) - k -1.
\end{align*}

 \item $e (X_t) = e (Y_1) + e (Y_2) - 2 e (Y_1 \cap Y_2) = e (Y_1) + e (Y_2) - 48$.
Then the result easily comes from
 $$e (X_t) = 2 (h^{1,1} (X_t) - h^{1,2} (X_t)),$$
  $$e (Y_i) = 2 (h^{1,1} (Y_i) - h^{1,2} (Y_i) +1).$$

\item Note that $\Pic (Y_i) \simeq H^2 (Y_i, \mbz)$, $\Pic (X_t) \simeq H^2 (X_t, \mbz)$ and $\Pic (D) \hookrightarrow H^2 (D, \mbz)$. Then it follows immediately from Corollary  \ref{two1}.
\item It also follows from Proposition \ref{6.4}.
\end{enumerate}
\end{proof}

We go over the examples of Calabi--Yau 3-folds which
are introduced in \cite{KaNa}, p. 408. They are constructed from two
copies of $\mbp^3$.  The examples  are all type of the above. One can see some phenomena that do not occur in the
$K3$ surface degeneration, which is the two-dimensional counterpart.

 Let $X_0 = Y_1 \cup Y_2$, where
$D = Y_1 \cap Y_2$ is a smooth quartic $K3$ surface in $\mbp^3$ and let $c$ and $c'$
be reduced divisors of $D$ which are composed of smooth curves. Let $c =  c_1 +
\cdots + c_s$ and $c' =  c'_1 + \cdots + c'_t$. We make $Y_1$ (resp, $Y_2$) by
blowing up successively along with the proper transforms of  $c_1, \cdots,
c_s$ (resp.  $c'_1, \cdots, c'_t$) in this order.
According to theorem 4.2 in \cite{KaNa},  there is a semistable degeneration of
a Calabi--Yau 3-fold which has $X_0$ as its central fiber if $c+c' \in |\mco_D (8)|$
and $X_0$ is K\"ahler. For example, if $s=1$ and $t=0$, then it is the `quick
example' in Section \ref{quick}. In \cite{KaNa}, they gave the topological Euler numbers.
We obtains more information by Corollary \ref{9.1}. 

The  Hodge numbers depend only on $c+c'$, not on individual $c$ or $c'$, or on the
ordering of blow-ups along their components. However $Y_1$  (resp.
$Y_2$) depends on $c$  (resp. $c'$) and the ordering of blow-ups.
Note that the ordering does not matter for case of surfaces   (\cite{Ku}).  So
different pairs of $Y_1$ and $Y_2$ may give Calabi--Yau 3-folds with the
same
Hodge numbers. One way to compare examples which have the same Hodge numbers is
to compare their Picard groups and Chern classes. Let us consider a pair of
examples. The first one is:
\cen{$c = c_1$ and $c'=c'_1$, where $c_1 \in |\mco_D  (5)|$ and  $c'_1
\in |\mco_D  (3)|$.}
Let $H_i$ be the divisor (on $Y_i)$ of the total transform of a hyperplane of
$\mathbb P^3$ and $E_1$ and $E'_1$ be the exceptional divisor over $c_1$ and $c'_1$  respectively.
Then
\begin{align*}
&\{  (l, l') \in \Pic  (Y_1) \times \Pic  (Y_2) \big| l|_D = l'_D\} = \\
&\{  (a H_1 + b E_1, cE'_1+  (a+5b-3c) H_2) \big| a, b, c \in \mbz\}.
\end{align*}
Let $\delta =   (-D, D) =   (E_1 - 4H_1,  4H_2 - E'_1)$, $e_1 =   (H_1, H_2)$ and $e_2
=   (5H_1 - E_1 , 0)$. Then
$$  (a H_1 + b E_1, cE'_1+  (a+5b-3c) H_2) =   (a+5b + c) e_1 -   (b+c) e_2 - c \delta.$$
By Corollary \ref{9.1}, the Picard group of $X_t$ is
$$\{  (l, l') \in \Pic  (Y_1) \times \Pic  (Y_2) \big| l|_D = l'_D\} / \delta = \langle e_1, e_2 \rangle  $$
up to torsion.
Then the cubic product is given by:
\cen{
$e_1^3 = 2$,
$e_1^2 e_2 = 5$,
$e_1 e_2^2 = 5$ and
$e_2^3 = 5$.
}
Let $M_1$ and $M'_1$ be a  fiber of the map $E_1 \ra c_1$ and $E'_1 \ra c'_1$ respectively. Then
$$H^4  (X_t, \mbz) = RG^4  (X_0, \mbz) = \langle   (H_1^ 2 - 4M_1, 0),   (M_1, M'_1)\rangle $$
up to torsion.
The bilinear form on $H^2  (X_t, \mbz) \times H^4  (X_t, \mbz)$ is given in the above basis by
$$\left  ( \begin{array}{ccc}
1& 0 \\
1& 1
\end{array}\right),$$
which is unimodular. So Condition \ref{consur}  for  this
example is satisfied. For the rest of examples, the unimodularities are all checked. We will not  give the explicit checking.

The products with the second Chern class are:
\cen{
$e_1 \cdot c_2  (X_t) = 32$ and
$e_2 \cdot c_2  (X_t) = 50$.}
The second smoothing which has the same Hodge numbers is prepared by setting:
\cen{$c = 0$ and $c'=c'_1 + c'_2 $, where $c'_1 \in |\mco_D  (5)|$ and
$c'_2 \in
|\mco_D  (3)|$}.
Since $c=0$, $Y_1 = \mbp^3$. \label{p1}
To get $ Y_2$, blow up $\mathbb P ^3$ along $c'_1$.  Let $
\tilde c'_2$ be the proper transform of $c'_2$.  Do the second blow-up
along $ \tilde c'_2$ to get $Y_2$. Let $F_1$ be the inverse image of $c'_1$.
Then the Picard group of $X_t$ is:
\cen{$\langle   (H_1, H_2),   (0, 5H_1 - F_1) \rangle $}
Let $f_1 =   (H_1, H_2)$ and $f_2 =   (0, 5H_1 - F_1)$, then the cubic
product is given by:
\cen{
$f_1^3 = 2$,
$f_1^2 f_2 = 5$,
$f_1 f_2^2 = 5$ and
$f_2^3 = 5$
}
and the products with the second Chern class are:
\cen{
$f_1 \cdot c_2  (X_t) = 32$ and
$f_2 \cdot c_2  (X_t) = 50$
}
Clearly, these two examples have the same cubic form and products with the Chern classes.
 Are they connected by smooth deformation? The answer is positive.
 Let  $X_0=Y_1 \cup Y_2$  be prepared
by blowing up $\mathbb P^3$ along $c$ (for $Y_1$) and $c' + l$ (for $Y_2$) where
$l$ is a smooth curve which is the last blow-up center for $Y_2$. Let $F$ be the
exceptional locus of $l$ in $Y_2$. Then $F$ is a $\mathbb P ^1$-bundle over $c$.
Let $X' \ra X$ be the blow-up along $F$, $E$ be the exceptional divisor, $\widetilde {Y_i}$ be the proper transform of $Y_i$ and $F' = E \cap \widetilde{Y_1}$. One can easily see that $E$ is a $\mbp^1$-bundle over $F'$. It can be shown by a result of \cite{Na} that $E$ can be smoothly contracted down to $F'$. Let
$X' \ra X''$ be the contraction.
 Then the final degeneration $X''$ has a
central normal crossing fiber such that one of its components  is the blow-up
of $\mathbb P^3$ along $c+l$ and the other one is the blow-up of $\mbp^3$ along
$c'$. So the center $c$ of the blow-up is moved from $Y_2$ to $Y_1$, keeping the total
space smooth and unchanged except for the central fiber. Let us remark this.

\begin{remark}\label{fido}
The top blow-up center is movable to the other side without affecting $X - X_0$.
\end{remark}

Now we know that the top blow-up curve can be moved without changing the generic
Calabi--Yau fibers.  But there are pairs of configurations of curves which
cannot be equal by a series of moving the top blow-up curves. What about those
degenerations? Do we have still the same thing? To answer this question, let us
consider another pair of examples.
The first one is  prepared by blowing up along
\cen{$c=0$ and $c' = h'_1+h'_2+h'_3$,}
where $h'_1 \in |\mco_D (5)| $, $h'_2 \in |\mco_D (2)| $ and $h'_3 \in |\mco_D (1)|$. Let
$F_1$, $F_2$ and $F_3$ be the inverse images of $h'_1$, $h'_2$ and $h'_3$ respectively.
Let $e_1 =  ( H_1, H_2)$,   $e_2
=  ( 0, 5 H_2 - F_1)$,  $e_3 =  ( 0, 2 H_2 - F_2)$. Then the Picard group of $X_t$
is:
\cen{$\langle e_1, e_2, e_3\rangle $}
and  the cubic form is given by:
$$\begin{array}{lllll}
\mu_{111}=2, &\mu_{112}= 5, & \mu_{113}=2,& \mu_{122}=5, & \mu_{123}=10,\\
 \mu_{133}=-4,& \mu_{222}=5,&\mu_{223}=10,&\mu_{233}=20,& \mu_{333}=-32,
\end{array}$$
where $\mu_{ijk} = e_i \cdot e_j \cdot e_k$.
The second one is prepared by  setting
\cen{$c=0$ and $c' = h'_1+h'_3+h'_2$}
Note that the order of blowing up  is changed and it cannot be changed to the
first one with any series of moving the top blow-up curve. Let $E_1$, $E_2$
and $E_3$ be the inverse images of $h'_1$, $h'_2$ and $h'_3$.
Let $f_1 =  ( H_1, H_2)$,   $f_2 =  ( 0, 5 H_2
- E_1)$,  $f_3 =  ( 0, 2 H_2 - E_2)$. Then the Picard group of $X_t$ is:
\cen{$\langle f_1, f_2, f_3\rangle $}
The cubic form is given by:
$$\begin{array}{lllll}
\nu_{111}=2, &\nu_{112}=5, &\nu_{113}=2, &\nu_{122}=5, &\nu_{123}=10, \\
\nu_{133}=-4, &\nu_{222}=5, &\nu_{223}=10, &\nu_{233}=20, &\nu_{333}=-40,
\end{array}$$
where $\nu_{ijk} = f_i \cdot f_j \cdot f_k$.
These two Calabi--Yau 3-folds have the same Hodge numbers $h^{1, 1} =3$ and $h^{1, 2}
=83$. Are their cubic forms also the same?  We can calculate the Aronhold S- and
T-invariants of these cubic forms to answer this question. For the first
one, the S- and T-invariants are:
\cen{$S_1 = 0$ and $T_1 = -86400$.}
For the second one:
\cen{$S_2 = 0$ and $T_2 = -38400$.}
Since the invariants are different, the cubic forms are different. The
Calabi--Yau 3-folds which are the smoothings of the above two normal crossings are not
diffeomorphic and accordingly can not be connected by smooth deformation. In
conclusion, the order of blow-ups matters. By changing the order of blow-ups,
we can easily produce Calabi--Yau 3-folds which have the same Hodge numbers but
have different cubic forms (so are non-diffeomorphic).
By calculating the cubic forms on the Picard groups, one can divide the examples of Kawamata and Namikawa into over 200 different deformation types.

We can also tell whether Calabi--Yau 3-folds, obtained by smoothing, are connected by projective flat deformation.
See the following theorem:

\begin{theorem}\label{hil} Let $Z_1$ and $Z_2$ be Calabi--Yau 3-folds. Then $Z_1$ and $Z_2$  belong
to the same Hilbert scheme of some projective space, and accordingly connected by
projective flat deformation, if and only if they have ample divisor $\rho_1$ and $\rho_2$ on $Z_1$ and $Z_2$ respectively such that
\cen{$\rho_1^3 = \rho_2^3$ and $\rho_1 \cdot c_2(Z_1) = \rho_2 \cdot c_2(Z_2).$}
\end{theorem}
\begin{proof}
Firstly we show the sufficiency.
By Kodaira vanishing theorem and Riemann--Roch formula, we have:
$$\dim(H^2( \mco_{Z_i}(n \rho_1))) = \chi(\mco_{Z_i}(n\rho_i)) = \frac{\rho_i^3}{6}n^3 + \frac{\rho_i \cdot c_2(Z_i)}{12} n$$
for any positive integer $n$. 
We know that $8 \rho_i$ is very ample on $Z_i$ (\cite{GaPu}). So $Z_1$ and $Z_2$  have embeddings to the same projective space $\mbp^N$ by the linear systems $|8\rho_1|$, $|8\rho_2|$ respectively and have the same Hilbert polynomials, where $N=\dim H^2( \mco_X(8\rho_1))-1=\dim H^2( \mco_X(8\rho_2))-1$. Therefore they belong to the same Hilbert Scheme. The second equivalent assertion follows from the fact that a Hilbert scheme is connected (\cite{Ha}). Conversely if $Z_1$ and $Z_2$ are connected by projective flat deformation, then they have the same Hilbert polynomial in some projective space. If we choose their hyperplane section $\rho_1$ and $\rho_2$ respectively, they satisfy the condition in the theorem.
\end{proof}
Since we can calculate the values required in the theorem, we can decide whether given two Calabi--Yau 3-folds, obtained by smoothing normal crossings are connected by projective flat deformation. In the next section, we will construct some new Calabi--Yau 3-folds and show how they are connected with other known Calabi--Yau 3-folds by projective flat deformation (Theorem \ref{100}).

\section{Calabi--Yau 3-folds with Picard number one and Fano 3-folds }\label{cyfano}

Among Calabi--Yau 3-folds, those with Picard number
one, that is for which the Picard lattice is generated by a single element,
bear special interest.
They have following numerical invariants:
\cen{$\rho^3$,  $\rho \cdot c_2 (Z)$ and $h^{1,2}(Z)$,}
where $Z$ is the Calabi--Yau and  $\rho$ is the unique ample generator of
$\Pic (Z)$. These three invariants determine the topological type of $Z$ if
it is simply-connected.

There is a close relation between Fano 3-folds and Calabi--Yau 3-folds with Picard number one.
In Section \ref{quick}, we saw that a normal crossing of  $\mbp^3$  and some blow-up of $\mbp^3$ is
smoothable to a Calabi--Yau 3-fold, which is a
degree 8 hypersurface in $\mbp (1,1,1,1,4)$:
$$X (8) \subset \mbp (1,1,1,1,4).$$
We can try the same construction on other pairs of Fano 3-folds.
Suppose that there are Fano 3-folds $V_1$ and $V_2$  such that
\begin{condition}\label{cyfanoc}
\begin{enumerate}
\item $V_1$ and $V_2$ contains copies of a $K3$ surface $D$ in their anticanonical systems.
\item There is a smooth curve $c$ in $|\mco_D({-}K_{V_1}{-}K_{V_2})|$.
\item  $\frac{{-}K_{V_1}}{r_1}|_D$ is linearly equivalent to $\frac{{-}K_{V_2}}{r_2}|_D$. We require this to guarantee K\"ahlerness of normal crossing varieties to be introduced later. \label{kkk}
\end{enumerate}
\end{condition}
Then we obtain a Calabi--Yau by smoothing $X_0$:
\begin{prop}\label{beforegood} Let $Y_1=V_1$ and $Y_2 \ra V_2$ be the blow along $c$. We make a normal crossing $X_0 = Y_1\cup_D Y_2$. Then it is smoothable to a Calabi--Yau 3-fold $X_t$.
Furthermore, $X_t$ has Picard number one if and only if 
$$\Min\{\rk\Pic(V_1), \rk\Pic(V_2)\}=1.$$
In that case, the Calabi--Yau 3-fold $X_t$ has
the following invariants:
\begin{enumerate}
\item $\rho^3 = (\frac{1}{r_1} + \frac{1}{r_2}) \delta$.
\item $\rho \cdot c_2 = -\frac{K_{V_1}}{r_1} \cdot c_2(V_1) - \frac{K_{V_2}}{r_2} \cdot c_2(V_2) + (r_1 + r_2) \delta $
\item $h^{1,2} = 22 + h^{1,2}({V_1}) + h^{1,2}(V_1) + \frac{1}{2} (r_1 + r_2)^2 \delta-\Max( \{h^2(V_1), h^2(V_2)\}$,
\end{enumerate}
where $\delta = -K_{V_1}^3 / r_1^2 {=} -K_{V_2}^3 / r_2^2$, $\rho$ is the ample generator of $\Pic(X_t)$ and $c_2=c_2(X_t)$.
\end{prop}
\begin{proof}D-semistability of $X_0$ follows from the fact that $c \in |\mco_D({-}K_{V_1}{-}K_{V_2})|$.
For the K\"ahlerness of $X_0$, we
need to find  ample divisors $H_i$ of $Y_i$ such that
$H_1|_D$ is linearly equivalent to  $H_2|_D$ as divisors of $D$. Let $\pi:Y_2 \ra V_2$ be the blow-up and
$E$  be the exceptional  divisor. Note that ${-}K_{Y_1}$ and ${-}K_{V_2}$ are  ample divisors
for $Y_1$ and $V_2$ respectively because they are Fano 3-folds.
For sufficiently large $n$, both of
\cen{$H_1 = - (n r_2-1)K_{Y_1}$ and
$H_2 = - (n r_1-1)\pi^* K_{V_2}  -  E$}
 are ample divisors for $Y_1$ and
$Y_2$ respectively. It is easy to check that $H_1|_D$ is linearly equivalent to
$H_2|_D$. One can easily show that the other conditions for smoothing are also satisfied.
So a Calabi--Yau 3-fold $X_t$ is obtained by smoothing $X_0$. Now we apply Corollary \ref{9.1} to calculate its invariants. We note that $\Pic(V_i)\simeq H^2(V_i, \mbz)$.
 WOLG, assume that $\rk \Pic(V_1) \leq \rk \Pic(V_2)$.
Note that
\begin{align*}
{\rm{rk}} ({\rm{im}} (H^2 (Y_1, \mbz) \oplus H^2 (Y_2, \mbz) \ra H^2 (D, \mbz))) =& \Max( \{\rk H^2(V_1, \mbz), \rk H^2(V_2, \mbz)\}\\
=&\rk H^2(V_2, \mbz)=h^2(V_2),
\end{align*}
where we used  Lefschetz hyperplane theorem.
So we have
$$h^{1,1}(X_t) =h^2(V_1)+(h^2(V_2)+1)-h^2(V_2)-1=h^2(V_1).$$
So $X_t$ has Picard number one if and only if 
$$1=h^2(V_1)=\rk \Pic(V_1)=\Min\{\rk\Pic(V_1), \rk\Pic(V_2)\}.$$
From now on we assume $\rk\Pic(X_t)=1$ and calculate the invariants.
 From (3) in Corollary \ref{9.1}, it is easy to show:
$\rho = (-\frac{K_{V_1}}{r_1},  -\pi^*\left(\frac{K_{Y_2}'}{r_2} \right))$. Then
$$
\rho^3 = \left(\frac{{-}K_{V_1}}{r_1}\right)^3  + \pi^*\left(-\frac{K_{V_2}}{r_2} \right)^3
       = \frac{{-}K_{V_1}^3}{r_1^3}  + \frac{{-}K_{V_2}^3}{r_2^3}
       = (\frac{1}{r_1} + \frac{1}{r_2}) \delta.
$$
Note that
\begin{align*}
\pi^*K_{Y_2} \cdot c_2(Y_2) &=K_{V_2} \cdot c_2(V_2) + K_{Y_2}\cdot c\\
                            &=K_{V_2} \cdot c_2(V_2) + K_{Y_2}\cdot ({-}K_{V_1}|_D{-}K_{V_2}|_D)\\
                            &=K_{V_2} \cdot c_2(V_2) - r_2(r_1 + r_2) \delta.
\end{align*}
So we have
\begin{align*}
\rho \cdot c_2(X_t) =& -\frac{K_{Y_1}}{r_1} \cdot c_2(Y_1) - \pi^*\left(\frac{K_{Y_2}}{r_2} \right) \cdot c_2(Y_2)\\
               =& -\frac{K_{V_1}}{r_1} \cdot c_2(V_1) - \frac{K_{V_2}}{r_2} \cdot c_2(V_2) + (r_1 + r_2) \delta.
\end{align*}
Finally we have,
\begin{align*}
h^{1,2}(X_t) &= 21+h^{1,2}(Y_1) + h^{1,2}(Y_2) -  h^2(V_2)\\
             &= 21+ h^{1,2}(V_1)+ h^{1,2}(V_2) + g(c) -h^2(V_2)\\
             &= 21 h^{1,2}(V_1)+ h^{1,2}(V_2) + (c^2 / 2+ 1) -\Max( \{h^2(V_1), h^2(V_2)\}\\
             &=h^{1,2}(V_1)+ h^{1,2}(V_2) + ({1}/{2}) (r_1 + r_2)^2 \delta -\Max( \{h^2(V_1), h^2(V_2)\}-1,
\end{align*}
where $g(c)$ is the genus of the curve $c$.

\end{proof}

So we can assume that $\rk\Pic(V_1)=1$. By Lefschetz hyperplane theorem, the restriction map
\begin{equation}\label{88}
\Pic(V_i) \ra \Pic(D)
\end{equation}
is primitive(\textit{i.e.} the cokernal has no torsion). So \eqref{kkk} in Condition \ref{cyfanoc} implies:
\begin{equation}\label{77}
-\frac{K_{V_1}^3}{r_1^2} {=} -\frac{K_{V_1}^3}{r_2^2}.
\end{equation}
Let $\mcv_1$ and $\mcv_2$ be Fano families that satisfy Equation \eqref{77}. It is not automatic that there are Fano 3-folds $V_1$, $V_2$ in $\mcv_1$, $\mcv_2$ respectively that satisfy Condition \ref{cyfanoc}. 

Let $V$ be a Fano 3-fold with Picard rank one from a Fano family $\mcv$. We can set $-\frac{K_{V}^3}{r^2}=2n-2$ for some integer $n \geq 0$. Let $\mck 3(n)$ be the moduli space of $K3$ surface with primitive ample divisor $h$ with $h^2=2n-2$ and $\mcp(\mcv)$ be the moduli space of all pairs $(V, D)$ for $V \in \mcv$ and a $K3$ surface $D \in |-K_V|$. Then we have a forgetful map:
\begin{equation}\label{88}
\varsigma: (V, D) \in \mbp(\mcv) \ra D \in \mck 3(n).
\end{equation}
We recall a result of A. Kovalev (\cite{Ko}):
\begin{lem}\label{kov} The image of the map $\varsigma$ is Zariski open (in particular, dense) in $\mck 3(n)$. When $\mcv$ is rigid, \textit{i.e} $h^{1,2}(V)=0$, the map is surjective.
\end{lem}
\begin{proof}
The first assertion is a special case of Theorem 6.45 in \cite{Ko}. For the second assertion, see the comment in Example 8.58 in \cite{Ko}.
\end{proof}

Firstly let us assume further that $\rk\Pic(V_2)=1$.
From the complete list of families of Fano 3-folds (see, for example, the appendix of \cite{IsPr}), we can find exactly 26 combinations of Fano families that satisfies Equation \eqref{77}. The above lemma implies that we can choose two Fano 3-folds $V_1$, $V_2$ from the two Fano families such that satisfy Condition \ref{cyfanoc} --- Note that $n_1=n_2$ and $V_1, V_2$ have Picard rank one.
We can calculate the invariants of the Calabi--Yau 3-folds  by Proposition \ref{beforegood}. It turns out that there is  a Calabi--Yau 3-fold in Table 1 of \cite{EnSt} whose invariants are overlapped.
Indeed there
is a \Q-Fano 4-fold $W$ for every such pair so that \begin{enumerate}
\item $V_1$ and $V_2$ are embedded into  $W$ and $V_1 + V_2 \in |-K_W|$.
\item There is a Calabi--Yau 3-fold $Z \subset W$ which is linearly
equivalent to $V_1 + V_2$.
\end{enumerate}
One can find a list for the correspondence between  Calabi--Yau 3-folds and \Q-Fano 4-folds in pages 83, 84 of \cite{Lee}.

Let us take an example.
The following two Fano families 
\begin{itemize}
\item $\mcv_1$: smooth hypersurfaces of degree 6 in  $\mbp(1^4,3)$,
\item $\mcv_2$: smooth hypersurfaces of degree 6 in $\mbp(1^3,2,3)$
\end{itemize}
satisfy Equation \eqref{77}, where $\mbp(1^4,3):=\mbp(1,1,1,1,3)$ and etc.

We can choose two Fano 3-folds, satisfying Condition \ref{cyfanoc}, from the above Fano families and calculate the invariants of the Calabi--Yau 3-fold, obtained by Proposition \ref{beforegood}:
\cen{$\rho^3 = 3, \rho \cdot c_2 = 42$ and $h^{1,2} = 103$.}
There is known a Calabi--Yau 3-fold with these invariants. It is a degree 6 hypersurface in a weighted projective space:
$X(6) \subset \mbp(1^4,2).$
Note that
$$X(6)(\subset \mbp(1^4,2)) = X(6,3)(\subset \mbp(1^4,2,3)),$$
which is linearly equivalent to
$$X(6,2)+X(6,1)(\subset \mbp(1^4,2,3))$$
as divisors in
$X(6)(\subset \mbp(1^4,2,3))$. Here $X(\alpha, \beta)$ denote a complete intersection of a degree $\alpha$ and a degree $\beta$ hypersurfaces.
So the Calabi--Yau, $X(6)(\subset \mbp(1^4,2)$, has a degeneration to a  normal crossing of $X(6,2)(\subset \mbp(1^4,2,3))$ and some blow-up of $X(6,1)(\subset \mbp(1^4,2,3))$.
Note that
$$X(6,2)(\subset \mbp(1^4,2,3)) = X(6)(\subset \mbp(1^4,3))$$
belongs to $\mcv_1$ and
$$X(6,1)(\subset \mbp(1^4,2,3)) = X(6)(\subset \mbp(1^3,2,3))$$
belongs to $\mcv_2$.

Now consider the case that $\rk\Pic(V_2) > 1$.
In this case, we can find around 30 combinations of Fano families that satisfies the above equation. Assuming that we can choose two Fano 3-folds, satisfying Condition \ref{cyfanoc}, from such a combination of Fano families, we can calculate the invariants of the Calabi--Yau 3-folds that may be obtained by Proposition \ref{beforegood} (one can find a partial list of it in page 93 of \cite{Lee}).

Comparing the invariants with those of known examples of Calabi--Yau 3-folds with Picard numbers (see, for example, Table 1 in \cite{EnSt}), we find that
there are no overlaps. So it is an interesting problem to check that there exist two Fano 3-folds, satisfying Condition \ref{cyfanoc}, in such a combination of Fano families. Let $\mcv_1$, $\mcv_2$ satisfy Equation \eqref{77}. Assume that $\mcv_1$ is rigid (\textit{i.e.} $h^{1,2}(V_1)=0$) so that we can use the second assertion of Lemma \ref{cov}. Then $\mcv_1= \{ V_1 \}$. Let $V_2 \in \mcv_2$ and $D$ be a $K3$ surface in $|-K_{V_2}|$. By Lemma, we have a copy of $D$ in $|-K_{V_1}|$ that satisfies Condition \ref{cyfanoc}. From this we have new examples of Calabi--Yau 3-fold with Picard number one.


\begin{example}[$\mcv_1$: $X_{22}$ (\S 12.2, Fano of index 1, degree 22 and Picard rank 1 in \cite{IsPr}) and $\mcv_2$: a Fano of index 1, degree 22 and Picard rank 2 (\S 12.3, No.\ 15 in \cite{IsPr})]\label{91}\

From these Fano 3-folds, we have a Calabi--Yau 3-fold with Picard number one by Proposition \ref{beforegood} with the following invariants:
\cen{$\rho^3 = 44, \rho \cdot c_2 = 92$ and $h^{1,2} = 68.$}
Denote this Calabi--Yau 3-fold by $\Xi_1$.
\end{example}


\begin{example}[$\mcv_1$: $X_{22}$ (\S 12.2, Fano of index 1, degree 22 and Picard rank 1 in \cite{IsPr}) and $\mcv_2$: a Fano of index 1, degree 22 and Picard rank 2 (\S 12.3, No.\ 16 in \cite{IsPr})]

The Calabi--Yau 3-fold with Picard number, obtained in this case, has the following invariants:
\cen{$\rho^3 = 44, \rho \cdot c_2 = 92$ and $h^{1,2} = 66.$}
Denote this Calabi--Yau 3-fold by $\Xi_2$.
\end{example}


\begin{example}[$\mcv_1$: $X_{22}$ (\S 12.2, Fano of index 1, degree 22 and Picard rank 1 in \cite{IsPr}) and $\mcv_2$: a Fano of index 1, degree 22 and Picard rank 3 (\S 12.4, No.\ 6 in \cite{IsPr})]

The Calabi--Yau 3-fold with Picard number, obtained in this case, has the following invariants:
\cen{$\rho^3 = 44, \rho \cdot c_2 = 92$ and $h^{1,2} = 64.$}

Denote this Calabi--Yau 3-fold by $\Xi_3$.
\end{example}



\begin{example}[$\mcv_1$: $X_{5}$ (\S 12.2, Fano of index 2, degree 5 and Picard rank 1 in \cite{IsPr}) and $\mcv_2$: a Fano of index 1, degree 10 and Picard rank 2 (\S 12.3, No.\ 4 in \cite{IsPr})]

The Calabi--Yau 3-fold with Picard number, obtained in this case, has the following invariants:
\cen{$\rho^3 = 15, \rho \cdot c_2 = 66$ and $h^{1,2} = 75.$}
Denote this Calabi--Yau 3-fold by $\Xi_4$.
\end{example}



\begin{example}[$\mcv_1$: A smooth quadric Fano 3-fold in $\mbp^4$ (\S 12.2, Fano of index 3, degree 2 and Picard rank 1 in \cite{IsPr}) and $\mcv_2$: a Fano of index 1, degree 6 and Picard rank 2 (\S 12.3, No.\ 2 in \cite{IsPr})]

The Calabi--Yau 3-fold with Picard number, obtained in this case, has the following invariants:
\cen{$\rho^3 = 8, \rho \cdot c_2 = 56$ and $h^{1,2} = 88.$}

Denote this Calabi--Yau 3-fold by $\Xi_5$.
\end{example}



\begin{example}[$\mcv_1$: A smooth quadric Fano 3-fold in $\mbp^4$ (\S 12.2, Fano of index 3, degree 2 and Picard rank 1 in \cite{IsPr}) and $\mcv_2$: $\mbp^1 \times S_1$ (\S 12.6,a Fano of index 1, degree 6 and Picard rank 10 in \cite{IsPr}(There is a typo in the list there))]

The Calabi--Yau 3-fold with Picard number, obtained in this case, has the following invariants:
\cen{$\rho^3 = 8, \rho \cdot c_2 = 56$ and $h^{1,2} = 60.$}

Denote this Calabi--Yau 3-fold by $\Xi_5$.
\end{example}



\begin{example}[ $\mcv_1$: $\mbp^3$ and $\mcv_2$: a Fano of index 1, degree 4 and Picard rank 2 (\S 12.3, No.\ 1 in \cite{IsPr})]

The Calabi--Yau 3-fold with Picard number, obtained in this case, has the following invariants:
\cen{$\rho^3 = 5, \rho \cdot c_2 = 50$ and $h^{1,2} = 92.$}

Denote this Calabi--Yau 3-fold by $\Xi_7$.
\end{example}


$\Xi_1, \cdots, \Xi_7$ are new examples of Calabi--Yau 3-folds with Picard number one that are topologically different from the known ones.  There are some hidden relation between these and some of known examples of Calabi--Yau 3-folds with Picard number one.
Let $Z_1$ be a quintic Calabi--Yau 3-fold in $\mbp^4$, 
 $Z_2$ be a Calabi--Yau 3-fold that is complete intersection of quadratic and quartic hypersurfaces in $\mbp^5$, 
 $Z_3$ be a Calabi--Yau 3-fold of $\rho^3=15, \rho \cdot c_2 = 66$, $h^{1,2}=76$ in Table 1 of \cite{EnSt} 
 and  $Z_4$ be a Calabi--Yau 3-fold of $\rho^3=44, \rho \cdot c_2 = 92$ , $h^{1,2}=65$ in Table 1 of \cite{EnSt}.
  Note that $Z_1$, $Z_2$, $Z_3$ and $Z_4$ have Picard number one.
\begin{theorem}\label{100} Consider some groups of Calabi--Yau 3-folds with Picard number one:
\cen{$\{\Xi_1, \Xi_2, \Xi_3, Z_4\} $, $\{\Xi_4, Z_3\} $, $\{\Xi_5, \Xi_6, Z_2\} $ and $\{ \Xi_7, Z_1 \}$.}
Then all the  Calabi--Yau 3-folds in each group belong to the same Hilbert scheme of some projective space and they are connected by (necessarily non-smooth) projective flat deformation.
\end{theorem}
\begin{proof}
 Note that the Calabi--Yau 3-folds in each group have the same values of
$$(\rho^3, \rho \cdot c_2 ).$$
So we are done by Theorem \ref{hil}.
\end{proof}


\end{document}